\documentclass[12pt,intlim,righttag]{article}
\usepackage{amsmath,amsthm}
\usepackage{amssymb}
\usepackage{amsfonts}

\newcommand{\la}{\langle}
\newcommand{\ra}{\rangle}

\newcommand\beq{\begin{equation}}
\newcommand\eeq{\end{equation}}

\theoremstyle{Theorem}

\theoremstyle{corollary}

\theoremstyle{remark}

\theoremstyle{definition}

\begin{document}
\title{New proofs of some results on BMO martingales using BSDEs}

\author{B. Chikvinidze $^{1),2)}$ and M. Mania $^{3),4)}$}

\date{~}
\maketitle

\begin{center}
$^{1)}$  Tbilisi State University, Chavchavadze Ave. 1, Tbilisi, Georgia\\

$^{2)}$  Institute of Cybernetics of Georgian Technical University,\\
         E-mail: beso.chiqvinidze@gmail.com \\

$^{3)}$  A. Razmadze Mathematical Institute of Tbilisi State University,\\

$^{4)}$  Georgian American University, Chavchavadze Ave. 17,Tbilisi, Georgia,\\
         E-mail: misha.mania@gmail.com
\end{center}

\begin{abstract}


{\bf Abstract.}
Using properties of backward stochastic differential equations
we give new proofs of some well known results on
BMO martingales  and improve some estimates of BMO norms.

\bigskip

 
\


\end{abstract}

\section{Introduction}
The BMO martingale theory is extensively used to study backward stochastic differential equations (BSDEs).
Some properties of BMO martingales was already used by Bismut[3]
when he discussed the existence and
uniqueness of a solution of some particular backward stochastic Riccati equations, choosing the BMO space
for the martingale part of the solution process. In the work of
Delbaen et al [5] conditions for the closedness of stochastic
integrals with respect to semimartingales in $L^2$ were established in relation to the problem of
hedging contingent claims and linear BSDEs. Most of this
conditions deal with BMO martingales and reverse H\"older
inequalities. BMO martingales naturally arise in
BSDEs with quadratic generators. When the generator of a BSDE has quadratic growth then the martingale part
of any bounded solution of the BSDE is a BMO martingale. This fact
was proved in [8, 11, 12, 13, 15, 16] under various degrees of
generality. Later, the BMO norms
were used to prove an existence, uniqueness and stability results
for BSDEs, among others in
[1, 2, 4, 11, 12, 16, 17].

The aim of this paper is to do the converse: to prove some
results on BMO martingales using the BSDE technique.

It is well known that if $M$ is a BMO martingale, then the
mapping $\phi \; : \; \mathcal{L}(P)\ni X \longrightarrow \tilde X=\langle{X,M}\rangle -X \in \mathcal{L}(\tilde P)$ is
an isomorphism of $BMO(P)$ onto $BMO(\tilde P)$, where $d\tilde P=\mathcal{E}_T (M)dP$. E. g.,
it was proved by Kazamaki [9, 10] that the inequality
$$
||\tilde X||_{BMO(\tilde P)}\leq C_K(\tilde M)\cdot ||X||_{BMO(P)}
$$
is valid for all $X\in BMO(P)$, where the constant $C_K(\tilde M)>0$ is
independent of $X$ but depends on the martingale $M$.
Using the properties of a suitable BSDE we prove this inequality
with a constant $C(\tilde M)$ which we express as a linear function of the
$BMO(\tilde P)$ norm of $\tilde M=\la M\ra - M$ and which is less than $C_K(\tilde M)$
for all values of this norm.

Using properties of BSDEs we prove also the well known equivalence
between BMO property, Muckenhoupt and reverse H\"older conditions
(Doleanse-Dade and Meyer [7], Kazamaki [10])
and obtain BMO norm estimates in terms of reverse H\"older and Muckenhaupt
constants.

\section{Reverse H\"older and Muckenhoupt conditions and relations with BSDEs}

We start with a probability space $\big( \Omega , \mathcal{F}, P \big)$, a finite time horizon \\
$0<T<\infty$ and a filtration
$F=(\mathcal{F}_t)_{0\leq t\leq T}$ satisfying the usual conditions of right-continuity and completeness.

We recall definitions of BMO martingales, Reverse H\"older and Muckenhaupt conditions
(see, e.g., Doleanse-Dade and Meyer [7], or Kazamaki [10]).

{\sc Definition 1.}
A continuous, uniformly integrable martingale
$(M_t, \mathcal{F}_t)$  with $M_0=0$ is said to be from the
class $BMO$ if
$$
||M||_{BMO}=\sup_\tau\Big\|E\big[\langle{M}\rangle _T-\langle{M}\rangle _\tau |\mathcal{F}_\tau \big]^{1/2}\Big\|_\infty <\infty,
$$
where the supremum is taken over all stopping times $\tau\in[0,T]$ and $\la M\ra$ is the square characteristic
of $M$.
\

Denote by $\mathcal E(M)$ the stochastic exponential of a
continuous local martingale $M$:
$$
\mathcal E_{t}(M)=\exp\{M_t-\frac{1}{2}\la M\ra_t\}.
$$

Throughout the paper we shall assume that $M$ is a continuous local martingale
with $\la M\ra_T<\infty$ $P$- a.s. This implies that
${\mathcal E}_T(M)>0$ $P$-a.s.  and let
 $\mathcal E_{\tau, T}(M)=\mathcal E_T(M)/\mathcal E_\tau(M)$.

\

{\sc Definition 2. }
Let $1<p<\infty $. $\mathcal{E}(M)$ is said to satisfy $(R_p)$ condition if the reverse H\"older inequality
$$
E\Big[ \big\{ \mathcal{E}_{\tau ,T}(M) \big\}^{p}  \Big| \mathcal{F}_\tau \Big]\leq C_p
$$
is valid for every stopping time $\tau $, with a constant $C_p>0$ depending only on $p$.

If $\mathcal{E}(M)$ is a uniformly integrable martingale then  by
the Jensen inequality we also have that
$E\Big[ \big\{ \mathcal{E}_{\tau ,T}(M) \big\}^{p}  \Big| \mathcal{F}_\tau \Big]\ge1$.

A condition dual to $(R_p)$ is the   Muckenhoupt condition $(A_p)$.

{\sc Definition 3. }
$\mathcal{E}(M)$ is said to satisfy  $(A_p)$ condition for $1<p<\infty $ if there is a constant $D_p>0 $
such that for every stopping time $\tau\in [0,T]$
$$
E\Big[ \big\{ \mathcal{E}_{\tau ,T}(M) \big\}^{-\frac {1}{p-1}}\Big| \mathcal{F}_\tau \Big]\leq D_p.
$$

Note that, since ${\mathcal E}(M)$ is a supermartingale, the
Jensen inequality implies the converse inequality

$$
E\Big[ \big\{ \mathcal{E}_{\tau ,T}(M) \big\}^{-\frac {1}{p-1}}  \Big| \mathcal{F}_\tau \Big]\geq
\Big\{ E \big[\mathcal{E}_{\tau ,T}(M)\big| \mathcal{F}_\tau \big] \Big\}^{-\frac {1}{p-1}}\geq 1.
$$

In this paper we shall consider only linear BSDEs of the type
$$
{Y_t =Y_0 - \int ^{t}_{0} [\alpha Y_{s} + \beta \psi _s ] d\langle{M}\rangle _s + \int ^{t}_{0} \psi _s d M_s + N_t}, \;\;\;
{Y_T =1},
$$
where $\alpha$ and $\beta$ are constants.
A solution of such BSDE we define as a triple $(Y,\psi,N)$,
with $\la N, M\ra=0$, from the space
$S^{\infty }\times BMO(P)\times H^2(P)$ equipped with the following norms
$$
||Y||_\infty=||Y_T^*||_{L^\infty},\;\;\;\text{where}\;\;\;\;
Y_T^*=\sup_{t\in[0,T]}|Y_t|,
$$
$$
||\psi\cdot M||_{BMO(P)} = \sup_\tau \Big\|E\Big[\int_\tau^T\psi_s^2d\langle{M}\rangle _s|\mathcal{F}_\tau \Big]^{1/2}\Big\|,
$$
$$
||N||_{H^2}=E^{\frac{1}{2}}[N]_T.
$$

Note that, since the martingale $M$ is assumed to be continuous,
only the latter term of this equation may have the jumps, i.e.,
$\Delta Y=\Delta N$.
In order to avoid the definition of BMO norms for right-continuous
martingales, we are using the $H^2$ norms for orthogonal martingale
parts. This is sufficient for our goals,
since the generators of equations under consideration does not depend on
orthogonal martingale parts.

Sometimes we call $Y$ alone the solution of BSDE, keeping in mind that
$\psi\cdot M+N$ is the martingale part of $Y$.

\

{\sc Lemma 1. }
Let $M$ be a continuous local martingale.
\\
{\sc \bf a) } $\mathcal{E}(M)$ satisfies $(R_p)$ if and only if there exists a bounded, positive solution of BSDE
\begin{equation}
\begin{cases}
{Y_t =Y_0 - \int ^{t}_{0} [\frac {p(p-1)}{2} Y_{s} + p \psi _s ] d\langle{M}\rangle _s + \int ^{t}_{0} \psi _s d M_s + N_t},\\

{Y_T =1}.
\end{cases}
\end{equation}
\\
{\sc \bf b) } $\mathcal{E}(M)$ satisfies $(A_p)$ if and only if there exists a bounded, positive solution of equation
\begin{equation}
\begin{cases}
{X_t =X_0 - \int ^{t}_{0} [\frac {p}{2(p-1)^2} X_{s} - \frac {1}{p-1} \varphi _s ] d\langle{M}\rangle _s + \int ^{t}_{0} \varphi _s d M_s + L_t},\\

{X_T =1}.
\end{cases}
\end{equation}

{\it Proof:}
{\sc \bf a) }
Let first show that if $\mathcal{E}(M)$ satisfies $(R_p)$ then the process
$Y_t=E\Big[ \big\{ \mathcal{E}_{t,T}(M) \big\}^p  \Big| \mathcal{F}_t \Big]$
is a solution of BSDE $(1)$. It is evident that $Y$ is a bounded positive process  and
that  $Y_t \big\{\mathcal{E} _t(M)\big\}^p$ is a uniformly integrable
 martingale. Therefore, since  $\mathcal{E}_t(M)>0$,
the process $Y$ will be a special semimartingale. Let
$Y_t = Y_0 + A_t + m_t$ be the canonical decomposition of $Y$, where
$m$ is a locally square integrable martingale and $A$
a predictable  process of bounded variation.
Using the Galtchouk-Kunita-Watanabe decomposition for $m$, we get
\begin{equation}\label{dec}
Y_t = Y_0 + A_t + \int ^{t}_{0} \psi _s d M_s + N_t,
\end{equation}
where $N$ is a local martingale orthogonal to $M$.

Now using the Ito formula we have
$$
Y_t \big\{\mathcal{E} _t(M)\big\}^p = Y_0 + \int _{0}^{t} \big[\frac {p(p-1)}{2}Y_{s}+p \psi _s\big]\big\{\mathcal{E} _s(M)\big\}^p d\langle{M}\rangle _s +
$$
\begin{equation}\label{dec2}
+\int _{0}^{t} \big\{\mathcal{E} _s(M)\big\}^p d A_s + \tilde m_t,
\end{equation}

where $\tilde m$ is a local martingale.

Because $Y_t \big\{\mathcal{E} _t(M)\big\}^p$ is a martingale, equalizing the part of  bounded variation to zero,
we obtain that

$$
A_t = - \int ^{t}_{0} \Big[\frac {p(p-1)}{2} Y_{s} + p \psi _s \Big] d\langle{M}\rangle _s,
$$
which implies that $Y_t=E\Big[ \big\{ \mathcal{E}_{t ,T}(M) \big\}^p  \Big| \mathcal{F}_t \Big]$ is a solution of equation $(1)$.

Now let equation $(1)$ admits a bounded positive solution $Y_t$.
Using the Ito formula for the process $Y_t \big\{\mathcal{E} _t(M)\big\}^p$ we get that
 $Y_t \big\{\mathcal{E} _t(M)\big\}^p$ is a local martingale. Hence it is a supermartingale, as a positive local martingale.
 Therefore, from the supermartingale inequality and the boundary condition $Y_T=1$ we obtain that
$E\Big[ \big\{ \mathcal{E}_{t,T}(M) \big\}^p \Big| \mathcal{F}_t \Big]\leq Y_t$. Because $Y$ is bounded,  this implies
 that $\mathcal{E} (M)$ satisfies $(R_p)$ condition.
\\
{\sc \bf b) }  The proof is similar to the proof of the part a), we only need to replace $p$ by $-\frac {1}{p-1}$.\qed

\

Let ${\mathcal E}(M)$ be a uniformly integrable martingale. Denote by $\tilde P$ a new probability
measure defined by $d\tilde P={\mathcal E}_T(M)dP$ and let $\tilde M= \la M\ra- M$.

Now we shall give a new proof of the well known equivalence
(Doleanse-Dade and Meyer [7], Kazamaki [10])
between BMO property, Muckenhoupt and reverse H\"older conditions.

\

{\sc \bf Theorem 1: }
Let $\mathcal{E} (M)$ be a uniformly integrable martingale. Then the following conditions are equivalent:
\\
{\sc \bf i). } $\tilde M \in BMO(\tilde P)$.
\\
{\sc \bf ii). } ${\mathcal E}(M)$ satisfies the $(R_p)$ condition
 for some $p>1 $.
\\
{\sc \bf iii). } $M \in BMO(P)$.
\\
{\sc \bf iv). } ${\mathcal E}(M)$ satisfies the $(A_p)$ condition for some $p>1 $.

\

{\it Proof:} For the sake of simplicity, in all proofs given here,
 we shall assume without loss of generality that all stochastic integrals are martingales,
otherwise one can use the localization arguments.

$i)\Longrightarrow ii)$
Let $\tilde M \in BMO(\tilde P)$.  According to Lemma 1 it
is sufficient to show that equation $(1)$ admits a bounded
positive solution for some $p>1$.
Let us rewrite  equation $(1)$ in terms of the $\tilde P$-martingale $\tilde{M}$:
$$
\begin{cases}
{Y_t =Y_0 - \int ^{t}_{0} [\frac {p(p-1)}{2} Y_{s} + (p-1) \psi _s ] d\langle{M}\rangle _s - \int ^{t}_{0} \psi _s d \tilde M_s + N_t},\\
{Y_T =1}.
\end{cases}
$$
Since $\la N,M\ra =0$, $N$ is a local $\tilde P$- martingale orthogonal to $\tilde M$.

Define the mapping $H \; : \; S^{\infty }\times {BMO(\tilde P)}\times H^2(\tilde P)$
into itself, which maps  $(y, \psi , n)\in S^{\infty }\times{BMO(\tilde P)}\times H^2(\tilde P)$
onto the solution $(Y,\Psi,N)$ of the BSDE (1), i.e.,
$$
Y_t=E^{\tilde P}\bigg[ 1+\int ^{T}_{t} \Big[\frac {p(p-1)}{2} y_{s} + (p-1)\psi _s \Big] d\langle{M}\rangle _s \bigg| \mathcal{F}_t \bigg]
$$
and
$$
-\int_0^t\Psi_s d \tilde M_s+N_t=E^{\tilde P}\bigg[ 1+\int ^{T}_{0} \Big[\frac {p(p-1)}{2} y_{s} + (p-1)\psi _s \Big] d\langle{M}\rangle _s \bigg| \mathcal{F}_t \bigg].
$$
We shall show that there exists $p>1$ such  that this mapping
is a contraction.

Let
$$
\delta Y=Y^1-Y^2, \; \delta y=y^1-y^2, \; \delta \Psi  =\Psi ^1-\Psi ^2, \; \delta \psi =\psi ^1-\psi ^2, \; \delta N=N^1-N^2.
$$
It is evident that $\delta Y_T =0$ and
$$\delta Y_t =\delta Y_0 - \int ^{t}_{0} \Big[\frac {p(p-1)}{2}\delta y_{s} + (p-1)\delta \psi _s \Big] d\langle{M}\rangle _s -
\int ^{t}_{0} \delta \Psi _s d \tilde M_s + \delta N_t.
$$
According to the Ito formula, applied for $(\delta Y_\tau )^2-(\delta Y_T)^2$
and taking conditional expectations we have
$$
(\delta Y_\tau )^2+E^{\tilde P}\bigg[\int ^{T}_{\tau } (\delta \Psi _s)^2 d\langle{M}\rangle _s \bigg| \mathcal{F}_\tau  \bigg] +
E^{\tilde P}\Big[ [\delta N]_T - [\delta N]_\tau  \Big| \mathcal{F}_\tau  \Big]=
$$
$$
= E^{\tilde P}\bigg[\int ^{T}_{\tau } p(p-1)\delta Y_{s} \delta y_{s} d\langle{M}\rangle _s \bigg| \mathcal{F}_\tau  \bigg]+
E^{\tilde P}\bigg[\int ^{T}_{\tau } 2(p-1)\delta Y_{s} \delta \psi _s d\langle{M}\rangle _s \bigg| \mathcal{F}_\tau  \bigg]
$$
and using elementary inequalities we obtain
$$
(\delta Y_\tau )^2+E^{\tilde P}\bigg[\int ^{T}_{\tau } (\delta \Psi _s)^2 d\langle{M}\rangle _s \bigg| \mathcal{F}_\tau  \bigg] +
E^{\tilde P}\Big[ [\delta N]_T - [\delta N]_\tau  \Big| \mathcal{F}_\tau \Big]\leq
$$
$$
\leq\frac {p(p-1)}{2} ||\tilde M||^2_{BMO(\tilde P)}\cdot ||\delta Y||^2_{\infty }+
\frac {p(p-1)}{2} ||\tilde M||^2_{BMO(\tilde P)}\cdot ||\delta y||^2_{\infty }+
$$
$$
+(p-1)||\tilde M||^2_{BMO(\tilde P)}\cdot ||\delta Y||^2_{\infty }+
(p-1)\Big\|\int \delta \psi d\tilde M\Big\|^2_{BMO(\tilde P)}.
$$
Because the right hand side of the inequality does not depend
on $\tau $, we will have

$$
\Big( 1 - p(p-1)||\tilde M||^2_{BMO(\tilde P)} - 2(p-1)||\tilde M||^2_{BMO(\tilde P)}\Big)||\delta Y||^2_{\infty }+
$$
$$
+\Big\|\int \delta \Psi d\tilde M\Big\|^2_{BMO(\tilde P)}+||\delta N||^2_{L^2(\tilde P)}\leq
$$
\begin{equation}\label{01}
\leq p(p-1)||\tilde M||^2_{BMO(\tilde P)}||\delta y||^2_{\infty } + 2(p-1)\Big\| \int \delta \psi d\tilde M\Big\|^2_{BMO(\tilde P)}.
\end{equation}
Since
$$
1 - (p-1)(p+2)||\tilde M||^2_{BMO(\tilde P)} <1
$$
for $p$ sufficiently close to $1$,  one can make the constant of $||\delta Y||^2_{\infty }$ in the left-hand side of (\ref{01}) positive
and  we finally  obtain the
inequality
$$
||\delta Y||^2_{\infty }+\Big\|\int \delta \Psi d\tilde M\Big\|^2_{BMO(\tilde P)}+||\delta N||^2_{H^2(\tilde P)}\leq
$$
\begin{equation}\label{lip}
\leq \alpha (p)\cdot ||\delta y||^2_{\infty }+\beta (p)\cdot \Big\|\int \delta \psi d\tilde M\Big\|^2_{BMO(\tilde P)},
\end{equation}

where
$$\alpha (p)=\frac {p(p-1)||\tilde M||^2_{BMO(\tilde P)}}{1-(p-1)(p+2)||\tilde M||^2_{BMO(\tilde P)}},
$$
$$
\beta (p)=\frac {2(p-1)}{1-(p-1)(p+2)||\tilde M||^2_{BMO(\tilde P)}}.
$$

It is easy to see that \; $\lim _{p\downarrow 1} \alpha (p)=\lim _{p\downarrow 1} \beta (p)=0$.
So, if we take $p^*$ such that $\alpha (p^*)<1$ and $\beta (p^*)<1$
we obtain that the mapping $H$ is a contraction
and there exists a unique solution $(Y,\Psi, N)$ of (1) in   $S^{\infty }\times  BMO(\tilde P)\times H^2(\tilde P)$.

Since $\alpha(p)$ and $\beta(p)$ are decreasing functions of
$p\in (1,\infty)$,  the norms $||Y||_\infty$ and
$||\Psi\cdot\tilde M||_{BMO(\tilde P)}$ are uniformly
bounded, as functions of $p$ for $p\in [1,p^*]$. Therefore,
for any $p\in [1,p^*]$ we have
\begin{equation}\label{ypos}
Y_t=E^{\tilde P}\bigg[ 1+\int ^{T}_{t} \Big[\frac {p(p-1)}{2} Y_{s} +
(p-1)\Psi _s \Big] d\langle{M}\rangle _s \bigg| \mathcal{F}_t \bigg]
\end{equation}
and
$$
Y_t\ge 1-\frac{p(p-1)}{2}||Y||_\infty ||\tilde M||_{BMO(\tilde P)}-\frac{p-1}{2}||\tilde M||_{BMO(\tilde P)}-
$$
$$
-\frac{p-1}{2}||\Psi \cdot \tilde M||_{BMO(\tilde P)}\ge0
$$
for some $p$ sufficiently close to $1$.
Hence, there exists a bounded positive solution of
equation (1) for some $p>1$, which implies that
 ${\mathcal E}(M)$ satisfies the $R_p$ condition, according to Lemma 1.

$ii)\Longrightarrow iii)$
Let $\mathcal{E} (M)$ be a uniformly integrable martingale and
satisfies the $(R_p)$ condition for some $p>1$.
Then the process $Y_t=E\Big[ \big\{ \mathcal{E}_{t,T}(M) \big\}^{p}  \Big| \mathcal{F}_t \Big]$
is a solution of equation (1) and satisfies the two-sided inequality
$$
1\le Y_t\le C_p.
$$

Using the Ito formula for $e^{-\beta Y_{T}}-e^{-\beta Y_{\tau }}$ and taking conditional
expectations we have
$$
e^{-\beta } - e^{-\beta Y_{\tau }} = \beta \frac {p(p-1)}{2} E\Big[\int ^{T}_{\tau } Y_{s} e^{-\beta Y_{s}} d\langle{M}\rangle _s \Big| \mathcal{F}_\tau  \Big]+
$$
$$
+E\Big[\int ^T_{\tau } e^{-\beta Y_{s}} \Big( \frac {\beta ^2}{2} \psi ^2_s + \beta p \psi _s \Big) d\langle{M}\rangle _s \Big| \mathcal{F}_\tau  \Big] +
\frac {\beta ^2}{2} E\Big[\int ^T_{\tau } e^{-\beta Y_{s}} d\langle{N^c}\rangle _s \Big| \mathcal{F}_\tau \Big]+
$$
$$
+E\Big[ \Sigma _{\tau < s \leq T} \big( e^{-\beta Y_s} - e^{-\beta Y_{s-}} + \beta e^{-\beta Y_{s-}}\Delta Y_s \big) \Big| \mathcal{F}_\tau  \Big].
$$
Since $\frac {\beta ^2}{2} \psi ^2_s + \beta p \psi _s \geq -\frac {p^2}{2}$,
$e^{-\beta Y_s} - e^{-\beta Y_{s-}} + \beta e^{-\beta Y_{s-}}\Delta Y_s \geq 0$ and $Y_t\ge1,$
taking $\beta >\frac {p}{p-1}$ we obtain the inequality

$$
\frac {p}{2}(\beta (p-1)-p)e^{-\beta C_p}E\Big[\langle{M}\rangle _{T}-\langle{M}\rangle _{\tau } \Big| \mathcal{F}_\tau  \Big] \leq e^{-\beta }-e^{-\beta C_p},
$$
which implies that
$$
||M||^2_{BMO(P)}\leq \frac {2(e^{\beta (C_p-1)}-1)}{p(\beta (p-1)-p)}
$$
for any $\beta>\frac{p}{p-1}.$

\

$iii)\Longrightarrow iv)$
If $M$ is a $BMO(P)$ martingale, then according to Lemma 1
it is sufficient to show that equation $(2)$ admits  bounded
positive solution for some $p>1$, which can be proved
similarly to the implication $i)\Longrightarrow ii)$.
By the same way one can show that for
the mapping $H$
$$
X_t=E\bigg[ 1+\int ^{T}_{t} \Big[\frac {p}{2(p-1)^2} x_{s} - \frac {1}{p-1} \varphi _s \Big] d\langle{M}\rangle _s \bigg| \mathcal{F}_t \bigg],
$$
where $-\int ^{t}_{0} \Phi _s d M_s + L_t$ is the martingale part of $X$,
the inequality (\ref{lip}) holds with
$$
\alpha (p)= \frac {p||M||^2_{BMO(P)}}{(p-1)^2-(3p-2)||M||^2_{BMO(P)}},
$$
$$
\beta (p)=\frac {2(p-1)}{(p-1)^2-(3p-2)||M||^2_{BMO(P)}},
$$
where \; $\lim _{p\rightarrow \infty } \alpha (p)=\lim _{p\rightarrow \infty } \beta (p)=0$.
So if we take $p$  large enough we obtain that the mapping $H$ is a contraction.

$iv)\Longrightarrow i)$
The proof is similar to the proof of  the implication $ii)\Longrightarrow iii)$.
In particular, for the BMO norm of $\tilde M$ the following
inequality holds
$$
||\tilde M||^2_{BMO(\tilde P)}\leq \frac {2(p-1)^2}{p(\beta -p)}\big( e^{\beta (D_p-1)}-1 \big)
$$
for any $\beta>p$, where $D_p$ is a constant from Definition 3.
\qed

\section{Girsanov's transformation of BMO martingales and BSDEs}

Let $M$ be a continuous local $P$-martingale such that $\mathcal{E}(M)$ is
a uniformly integrable martingale and let
$d\tilde P=\mathcal{E}_T (M)d P$.
To each continuous local martingale $X$ we associate the process
 $\tilde X=\la X,M\ra-X$, which is a local $\tilde P$-martingale
 according to Girsanov's theorem. We denote this map by
 $\varphi:\mathcal{L}(P)\to\mathcal{L}(\tilde P)$,
where
$\mathcal{L}(P)$ and $\mathcal{L}(\tilde P)$ are classes of $P$ and
$\tilde P$ local martingales.

Let consider the process
\begin{equation}\label{y}
Y_t=E^{\tilde P}\big[\langle{X}\rangle _{T}-\langle{X}\rangle _{t} \big| \mathcal{F}_t \big]
=E\big[{\mathcal E}_{t,T}(M)(\langle{X}\rangle _{T}-\langle{X}\rangle _{t}) \big| \mathcal{F}_t \big].
\end{equation}
Since $\la\tilde X\ra=\la X\ra$ under either probability measure, it is evident that
$$
||Y||_{\infty }=||\tilde X||^2_{BMO(\tilde P)}.
$$
Let $M \in BMO(P)$.
According to Theorem 1  condition $(R_p)$ is satisfied
 for some $p>1$.  The $(R_p)$
condition and conditional energy inequality (Kazamaki [10], page 29)
imply that for any $X \in BMO(P)$ the process $Y$ is
bounded, i.e., $\varphi$ maps $BMO(P)$ into $BMO(\tilde P)$.
Moreover, as proved by Kazamaki [9, 10], $BMO(P)$ and $BMO(\tilde P)$
are isomorphic under the mapping $\phi$
and  for all $X\in BMO(P)$ the inequality
\begin{equation}\label{kaza1}
||\tilde X||^2_{BMO(\tilde P)} \leq C^2_K(\tilde M)\cdot||X||^2_{BMO(P)}
\end{equation}
is valid, where
\begin{equation}\label{ck}
C^2_K(\tilde M)= 2p \cdot 2^{1/p} \sup _{\tau }\Big\| E^{\tilde{P}}\Big[ \big\{ \mathcal{E}_{\tau ,T}(\tilde{M}) \big\}^{-\frac {1}{p-1}}\Big| \mathcal{F}_\tau \Big] \Big\|_{\infty } ^{(p-1)/p},
\end{equation}
and $p$ is such that
\begin{equation}\label{ck2}
||\tilde M||_{BMO(\tilde P)}<\sqrt{2}(\sqrt{p}-1).
\end{equation}
Note that  the similar inequality holds for the inverse mapping $\phi^{-1}$.

Now we give an alternative proof of this assertion, which improves also the constant
in the inequality (\ref{kaza1}).

\

{\sc \bf Theorem 2. }
If $M\in BMO(P)$, then $\phi \; : \; X\rightarrow \tilde X$ is an isomorphism of $BMO(P)$ onto $BMO(\tilde P)$.
In particular, the inequality
$$
\frac{1}{\Big( 1 + \frac{\sqrt 2}{2} ||M||_{BMO(P)} \Big)} ||X||_{BMO(P)} \leq ||\tilde X||_{BMO(\tilde P)} \leq
$$
\begin{equation}\label{main}
\leq \Big( 1 + \frac{\sqrt 2}{2} ||\tilde M||_{BMO(\tilde P)} \Big)||X||_{BMO(P)}.
\end{equation}
is valid for any $X\in BMO(P)$.

\

{\it Proof:}
Similarly to Lemma 1 one can show that for any $X\in BMO(P)$
the process
$Y$(defined by (\ref{y})) is a positive bounded solution of
the $BSDE$
\begin{equation}
\begin{cases}
{Y_t=Y_0-\langle{X}\rangle _{t}-\int ^{t}_{0}\varphi _s d\langle{M}\rangle _{s}+\int ^{t}_{0}\varphi _s d M_s + L_t},\\
{Y_T=0}.
\end{cases}
\end{equation}

Applying the Ito formula for $(Y_\tau +\varepsilon )^p-(Y_T+\varepsilon )^p$
where $0<p<1, \; \varepsilon >0$  and taking conditional expectations we obtain
$$
\big( Y_\tau  + \varepsilon  \big)^p - \varepsilon ^p = E\Big[\int ^{T}_{\tau } p(Y_s+\varepsilon )^{p-1}d\langle{X}\rangle _{s} \Big| \mathcal{F}_\tau  \Big] + \frac {p(1-p)}{2} E\Big[\int ^{T}_{\tau } (Y_s+\varepsilon )^{p-2} d \langle{L^c}\rangle _{s} \Big| \mathcal{F}_\tau  \Big]+$$
$$+E\Big[\int ^{T}_{\tau } \Big( \frac {p(1-p)}{2} (Y_s+\varepsilon )^{p-2}\varphi ^2_s +p(Y_s+\varepsilon )^{p-1} \varphi _s \Big) d\langle{M}\rangle _{s} \Big| \mathcal{F}_\tau  \Big] - $$
\begin{equation}\label{eps}
-E\Big[ \Sigma _{\tau <s\leq T} \big( (Y_s+\varepsilon)^p - (Y_{s-}+\varepsilon)^p - p (Y_{s-}+\varepsilon)^{p-1}\Delta Y_s \big) \Big| \mathcal{F}_\tau  \Big].
\end{equation}

Because  $f(x) = x^p$  is concave for $p\in(0,1)$,
the last term in (\ref{eps}) is positive. Therefore, using
the inequality
$$
\frac {p(1-p)}{2} (Y_s+\varepsilon )^{p-2}\varphi ^2_s +
p(Y_s+\varepsilon )^{p-1} \varphi _s+\frac{p}{2(1-p)}(Y_s+\varepsilon)^p\ge0
$$
from (\ref{eps}) we obtain
$$
(Y_\tau + \varepsilon )^p-\varepsilon^p\geq E\Big[\int ^{T}_{\tau } p(Y_s+\varepsilon )^{p-1}d\langle{X}\rangle _{s} \Big| \mathcal{F}_\tau  \Big]-
$$
\begin{equation}\label{eps2}
- \frac {p}{2(1-p)}E\Big[\int ^{T}_{\tau } (Y_s+\varepsilon )^p d\langle{M}\rangle _{s} \Big| \mathcal{F}_\tau  \Big].
\end{equation}

Since $0<p<1$
$$
p \big( ||Y||_{\infty } +\varepsilon\big)^{p-1}E\Big[\la X\ra_T-\la X\ra_\tau \Big| \mathcal{F}_\tau \Big]\le E\Big[\int ^{T}_{\tau } p(Y_s+\varepsilon )^{p-1}d\langle{X}\rangle _{s} \Big| \mathcal{F}_\tau \Big],
$$
from (\ref{eps2}) we have
$$
p \big( ||Y||_{\infty } +\varepsilon \big)^{p-1}
E\Big[\la X\ra_T-\la X\ra_\tau \Big| \mathcal{F}_\tau \Big]\leq
\big( Y_\tau + \varepsilon  \big)^p -\varepsilon^p + \frac {p}{2(1-p)}E\Big[\int ^{T}_{\tau } (Y_s+\varepsilon )^p d\langle{M}\rangle _{s} \Big| \mathcal{F}_\tau  \Big]
$$
and taking norms in the both sides of the latter inequality we obtain
$$
p \big( ||Y||_{\infty } +\varepsilon \big)^{p-1}\cdot ||X||^2_{BMO(P)}\leq \big( ||Y||_{\infty }
 +\varepsilon \big)^p-\varepsilon^p + \frac {p}{2(1-p)}\big( ||Y||_{\infty } +\varepsilon \big)^p \cdot  ||M||^2_{BMO(P)}.
$$

Taking  the limit when $\varepsilon \rightarrow 0$ we will have that for all $p\in (0,1)$
$$
||X||^2_{BMO(P)}\leq \Big( \frac {1}{p} + \frac {1}{2(1-p)}||M||^2_{BMO(P)} \Big)\cdot ||Y||_{\infty }.
$$
Therefore,
$$
||X||^2_{BMO(P)}\leq \min_{p\in (0,1)}\Big( \frac {1}{p} + \frac {1}{2(1-p)}||M||^2_{BMO(P)} \Big)\cdot ||Y||_{\infty }=
$$
\begin{equation}\label{min}
=\Big( 1 + \frac {\sqrt 2}{2} ||M||_{BMO(\tilde P)} \Big)^2\cdot ||Y||_{\infty },
\end{equation}

since the minimum of the function
$f(p)=\frac {1}{p} + \frac {1}{2(1-p)}||M||^2_{BMO(\tilde P)}$
is attained for $p^*=\sqrt 2 / (\sqrt 2 + ||\tilde M||_{BMO(\tilde P)})$ and
$f(p^*)=\Big( 1 + \frac {\sqrt 2}{2} ||M||_{BMO(\tilde P)} \Big)^2$.

Thus, from (\ref{min})
$$
\frac{1}{\Big( 1 + \frac {\sqrt 2}{2} ||M||_{BMO(P)} \Big)}||X||_{BMO(P)}\le ||\tilde X||_{BMO(\tilde P)}.
$$

Now we can use inequality (\ref{min}) for the Girsanov transform of $\tilde X$.

Since $dP/d\tilde P={\mathcal E}^{-1}_T(M)={\mathcal E}_T(\tilde M)dP$, $\tilde M, \tilde X\in BMO(\tilde P)$ and
$$
\varphi(\tilde X)=\tilde X -\la\tilde X,\tilde M\ra = X,
$$
from (\ref{min}) we get the inverse inequality:
\begin{equation}\label{ineq}
||\tilde X||_{BMO(\tilde P)}\leq \Big( 1 + \frac {\sqrt 2}{2} ||\tilde M||_{BMO(\tilde P)} \Big)||X||_{BMO(P)}.
\end{equation}
\qed

Let us compare the constant
$$
C(\tilde M)=  1 + \frac {\sqrt 2}{2} ||\tilde M||_{BMO(\tilde P)}
$$
from (\ref{main}) with the corresponding constant $C_K(\tilde M)$
from (\ref{kaza1}) (Kazamaki [10]).

Since
$$
E^{\tilde{P}}\Big[ \big\{ \mathcal{E}_{\tau ,T}(\tilde{M}) \big\}^{-\frac {1}{p-1}}\Big| \mathcal{F}_\tau \Big]\geq 1,
$$
the constant $C_K(\tilde M)$ is more than $\sqrt{2p}$,
where
$p$ is such that \\
 $||\tilde M||_{BMO(\tilde P)}<\sqrt{2}(\sqrt{p}-1)$.
Since the last inequality is equivalent to
the inequality $p>\big(1+\frac {\sqrt {2}}{2}||\tilde M||_{BMO(\tilde P)}\big)^2$,
we obtain that at least
$$
C^2(\tilde M)\le \frac{1}{2}C^2_K(\tilde M).
$$

From inequality (\ref{main}) it follows the following simple corollary, which can not be
deduced from inequality (\ref{kaza1}).

{\bf Corollary.} Let $(M^n, n\ge1)$ be a sequence of $BMO(P)$ martingales such that
$\lim_{n\to\infty}||M^n||_{BMO(P)}=0$. Let $dP^n={\mathcal E}_T(M^n)dP$ and ${\tilde X}^n=X-\la X, M^n\ra$.
Then for any $X\in BMO(P)$
$$
\lim_{n\to\infty}||\tilde X^n||_{BMO(P^n)} = ||X||_{BMO(P)}.
$$
{\it Proof.}  The second inequality of (\ref{main}), applied for $X=M^n$ and $M=M^n$ gives
$$
||\tilde M^n||_{BMO(P^n)}
\leq \Big( 1 + \frac{\sqrt 2}{2} ||\tilde M^n||_{BMO(P^n)} \Big)
||M^n||_{BMO(P)}.
$$
Therefore,
$$
\frac{1}{ \frac{\sqrt 2}{2} + 1/||\tilde M^n||_{BMO(P^n)}} \leq ||M^n||_{BMO(P)},
$$
which implies that   $\lim_{n\to\infty}||\tilde M^n||_{BMO(P^n)}=0$.
Now, passing to the limit in the two-sided inequality (\ref{main}) we obtain
$$
 ||X||_{BMO(P)}\leq\lim_{n\to\infty}||\tilde X^n||_{BMO(P^n)}
 \leq ||X||_{BMO(P)}.\qed
$$
{\bf Remark.} Note that the converse of Theorem 2 is also true. I.e., if $M$ is a continuous local martingale and ${\mathcal E}(M)$
is a uniformly integrable martingale, Schachermayer [18] proved that if $M\notin BMO(P)$ then the map $\varphi$ is not an isomorphism
from $BMO(P)$ into $BMO(\tilde P)$.

\newpage

\begin{center}
REFERENCES
\end{center}

\
\\
{\sc \bf [1] } S. Ankirchner, P. Imkeller, and G. Reis, " Classical and
variational differentiability of bsdes with quadratic growth"
{\it Electronic Journal of Probability}, Vol. 12 (2007),
pp. 1418-1453.

\
\\
{\sc \bf [2] } P. Barrieu, N. Cazanave and  N. El Karoui, "Closedness results for BMO semi-martingales and application to quadratic BSDEs,"
               { \it Comptes Rendus Mathematique}, Vol. 346, pp. 881-886, 2008.

\
\\
{\sc \bf [3] } J. M. Bismut, "Controle des systemes lineaires quadratiques: applications de l'integrale stochastique,"
               { \it Seminaire de Probabilites XII (eds.: C. Dellacherie, P. A. Meyer, and M. Weil), Lecture Notes in Mathematics 649},
               Springer-Verlag, Berlin/Heidelberg, pp. 180-264, 1978.

\
\\
{\sc \bf [4] } B. Chikvinidze, "Backward stochastic differential equations with a convex generator,"
               { \it Georgian Mathematical Journal}, Vol. 19, pp. 63-92, 2012.

\
\\
{\sc \bf [5] } F. Delbaen, P. Monat, W. Schachermayer, M. Schweizer and C. Stricker, "Weighted norm inequalities and hedging in incomplete markets,"
                { \it Finance and Stochastics }, Vol. 1, pp. 181-227, 1997.

\
\\
{\sc \bf [6] } F. Delbaen and S. Tang, "Harmonic analysis of stochastic equations and backward stochastic differential equations,"
                { \it Probability Theory and Related Fields }, Vol. 146, pp. 291-336, 2010.

\
\\
{\sc \bf [7] } C. Doleans-Dade and P. A. Meyer, "Inegalites de normes avec poids,"
                { \it Universite de Strasbourg Seminaire de Probabilites }, XIII,  pp. 313-331, 1979.

\
\\
{\sc \bf [8] } Y. Hu, P. Imkeller, and M. M\"uller, "Utility maximization in incomplete markets,"
               { \it Annals of Applied Probability }, Vol. 15, pp. 1691-1712, 2005.
 
\
\\
{\sc \bf [9] } N. Kazamaki, "On transforming the class of BMO-martingales by a change of law,"
                 { \it Tohoku Mathematical Journal }, Vol. 31, pp. 117-125, 1979.

\
\\
{\sc \bf [10] } N. Kazamaki, { \it Continuous Exponential Martingales and BMO }, vol. 1579 of
                { \it Lecture Notes in Mathematics }, Springer, Berlin-Heidelberg, 1994.

\
\\
{\sc \bf [11] } M. Kohlmann and S. Tang, "Minimization of risk and linear quadratic optimal control theory,"
                { \it SIAM Journal on Control and Optimization }, Vol. 42, pp. 1118-1142, 2003.

\
\\
{\sc \bf [12] } M. Mania and R. Tevzadze, "A semimartingale Bellman equation and the variance-optimal martingale measure,"
                { \it Georgian Mathematical Journal }, Vol. 7, pp. 765-792, 2000.

\
\\
{\sc \bf [13] } M. Mania and  M. Schweizer, "Dynamic exponential indifference valuation,"
                 { \it Annals of Applied Probability }, Vol. 15, pp. 2113-2143, 2005.

\
\\
{\sc \bf [14] } M. Mania and R. Tevzadze, "Martingale equation of exponential type,"
                { \it Electronic communication in probability }, Vol. 11, pp. 206-216, 2006.

\
\\
{\sc \bf [15] } M. Mania, M. Santacroce and R. Tevzadze,  "A semimartingale BSDE related to the
minimal entropy martingale measure",{\it Finance and Stochastics}, Vol. 7,  No. 3, pp. 385-402, 2003.

\
\\
{\sc \bf [16] } R. Tevzadze, "Solvability of backward stochastic differential equations with quadratic growth,"
                { \it Stochastic Processes and their Applications }, Vol. 118, pp. 503–515, 2008.

\
\\
{\sc \bf [17] } M.A. Morlais, "Quadratic BSDEs driven by a continuous martingale and application to
utility maximization problem", {\it Finance and Stochastics}, Vol. 13, No. 1, pp. 121-150, 2009.

\
\\
{\sc \bf [18] } W. Schachermayer, "A characterization of the closure of $H_{\infty }$ in BMO,"
                { \it Seminaire de Probabilites XXX, Lecture Notes in Mathematics 1626}, Springer, Berlin, pp. 344-356, 1996.

\end{document}